\documentclass[12pt]{amsart}
\usepackage{amsmath}
\usepackage{amssymb}
\usepackage{mathrsfs}
\usepackage[a4paper]{geometry}
\theoremstyle{plain}
\newtheorem{teo}{Theorem}
\newtheorem{pro}{Proposition}
\newtheorem{lem}{Lemma}

\theoremstyle{definition}

\theoremstyle{remark}

\newtheorem*{nota}{Remark}
\newtheorem*{note}{Remarks}
\newcommand{\la}{\langle}
\newcommand{\ra}{\rangle}

\newcommand{\p}{\mathbb{P}}
\newcommand{\re}{\mathbb{R}}

\date{\today}
\title{Radial Dunkl Processes: Existence, uniqueness and hitting time} 
\begin{document}
\maketitle
\centerline{NIZAR DEMNI\footnote{IRMAR, Rennes 1 University; email: nizar.demni@univ-rennes1.fr.}}  

\begin{abstract} 
We give shorter proofs of the following known results: the radial Dunkl process associated with a reduced system and a strictly positive multiplicity function is the unique strong solution for all times $t$ of a stochastic differential equation with a singular drift (see \cite{Scha} for the original proof via Dirichlet forms), the first hitting time of the Weyl chamber by a radial Dunkl process is finite almost surely for small values of the multiplicity function. Compared to the original proofs, ours give more information on the behavior of the process. More precisely, the first proof allows to give a positive answer to a conjecture announced by Gallardo and Yor in \cite{Chy} for the finiteness of the total length of the jumps performed by the Dunkl process during a finite amount of time, while the second one shows that the process hits almost surely the wall corresponding to the simple root with a small multiplicity value. Further results on the jumps of the Dunkl process (communicated by D. L\'epingle) and and additional references are included.    
\end{abstract}

\section{Preliminaries and examples}
To be self-contained, we begin by pointing out some facts on root systems and radial Dunkl processes. The reader is referred to \cite{Chy} and \cite{Hum} for more details. 
Let $(V, \la,\ra)$ be a finite-dimensional Euclidean space of dimension $m$. A \emph{reduced} root system $R$ is a finite set of non zero vectors in $V$ such that : 
\begin{itemize}
\item[1] $R \cap \re \alpha = \{\alpha,-\alpha\}$ for all $\alpha \in R$, \\
\item[2] $\sigma_{\alpha}(R) = R$, 
\end{itemize} where $\sigma_{\alpha}$ is the reflection with respect to the hyperplane $H_{\alpha}$ orthogonal to $\alpha$: 
\begin{equation*}
\sigma_{\alpha} (x) =  x - 2\frac{\la\alpha, x \ra}{|\alpha|^2} \alpha,\, |\alpha|^2 := \la \alpha,\alpha\ra \quad x \in V.
 \end{equation*}
A simple system $S$ is  a basis of $\textrm{span}(R)$ which induces a total ordering in $R$. A  root $\alpha$ is positive if it is a positive linear combination of elements of $S$. The set of positive roots is called a positive system and is denoted by $R_+$. The (finite) reflection group $W$ is the group generated by all the reflections $\sigma_{\alpha}$ for $\alpha \in R$ and acts on $R$ via the relation (\cite{Hum})
\begin{equation*}
\sigma_{\alpha}\sigma_{\eta}\sigma_{\alpha} = \sigma_{\sigma_{\alpha}(\eta)}, \, \alpha,\eta \in R.
\end{equation*} 
Given a root system $R$  with positive and simple systems $R_+,S$,  define the {\it positive Weyl chamber} $C$ by: 
\begin{equation*}
C := \{x \in V, \, \langle \alpha , x \rangle > 0 \, \forall \, \alpha \in R_+\} = \{x \in V, \, \langle \alpha , x \rangle > 0 \, \forall \, \alpha \in S\} 
 \end{equation*} 
 and $\overline{C},\partial C$ its closure and boundary respectively. The radial Dunkl process $X^W$ is defined as the $\overline{C}$-valued continuous-paths Markov process whose generator is given by : 
\begin{equation*}
\mathscr{L}_k^Wu(x) = \frac{1}{2}\Delta u(x) + \sum_{\alpha \in R_+}k(\alpha)\frac{\langle \alpha,\nabla u(x)\rangle}{\langle\alpha, x\rangle} 
\end{equation*} 
where $u \in C^2(\overline{C})$ satisfies the boundary conditions $\langle\nabla u(x), \alpha\rangle = 0$ for all $x \in H_{\alpha},\, \alpha \in R_+$, and $k(\alpha) \geq  0 $ is a multiplicity function (a $W$-invariant function).  

In order to motivate the reader and prepare for the first result, we exhibit some known examples. The first and easiest one corresponds to $V = \re, R = B_1 = \{\pm 1\}$. There is only one orbit so that $k(\alpha) := k \geq 0$ and $X^W$ is a \emph{Bessel} process (\cite{Rev}) of \emph{index} $\nu = k - 1/2$. When $k > 0$ and $X_0^W= x \geq 0$, it is the unique strong solution of the following stochastic differential equation with singular drift: 
\begin{equation*}
dX^W_t = dB_t + \frac{k}{X_t^W}dt, \quad t \geq 0,
\end{equation*}
where $B$ is a standard Brwonian motion. When $k=0$, $X^W$ is a reflected Brownian motion. A multivariate known example is provided by the so-called $A$-type root systems:
\begin{equation*}
R = A_{m-1} = \{\pm (e_i - e_j), \, 1 \leq i < j \leq m-1\}, 
\end{equation*} 
with positive and simple systems given by : 
\begin{equation*}
R_+ = \{e_i - e_j, \,1\leq i < j \leq m \}, \quad S = \{e_i - e_{i+1}, \,1\leq i  \leq m-1 \},
\end{equation*}
where $(e_i)_{1 \leq i \leq m}$ is the canonical basis of $\re^m$. In this case, $V = \re^m$, the span of $R$ is the hyperplane of $\re^m$ consisting of vectors whose coordinates sum to zero and $C = \{x \in \re^m, x_1 > \dots > x_m\}$. Besides, there is only one orbit so that $k(\alpha) := k \geq 0$ and $X^W = (X^{W,i})_{1 \leq i \leq m}$ satisfies :  
\begin{equation} 
\label{Khil}
dX_t^{W,i} = d\nu_t^i + k \sum_{j \neq i}\frac{dt}{X_t^{W,i} - X_t^{W,j}} \quad 1 \leq i \leq m , \quad t < \tau
\end{equation}
with $X_0^{W,1} > \dots > X_0^{W,m}$, where $(\nu^i)_i$ are independent BMs and $\tau$, \emph{the first collision time}, is defined by
\begin{equation*}
\tau := \inf\{t, \, X_t^{W,i} = X_t^{W,j} \, \textrm{for some} \, (i,j)\}.
\end{equation*}
This process fits the eigenvalues of symmetric and Hermtian Brownian motions when $k=1/2$ and $1$ respectively (\cite{Dyson}). For general values $k > 0$, it was deeply studied in \cite{Cepa} (see also \cite{Cepa1}) where it was proved that (\ref{Khil}) has a unique strong solution for all $t \geq 0$ any $X_0^{W} \in \overline{C}$. When $k=0$, $X^W$ is a reflected Brownian motion in the Weyl chamber of type $A$ and we refer the reader to \cite{DL} for properties of the reflected process for general root systems. Another, yet less familiar, instance of $X^W$ is a generalization to an arbitrary set of parameters of singular values of Wishart and Laguerre processes. 
To see this, recall from \cite{Bru},\cite{Dem} that their eigenvalues processes $(\lambda_i(t), 1 \leq i \leq m)_{t \geq 0}$ are unique strong solutions of: 
\begin{equation}
\label{Lag}
d\lambda_i(t) = 2\sqrt{\lambda_i(t)}\, d\nu_i(t) + \beta\left[ \delta + \sum_{k \neq i}\frac{\lambda_i(t) + \lambda_k(t)}{\lambda_i(t) - \lambda_k(t)}\right]dt, \quad  t < \tau \wedge R_0,
\end{equation}
where $\beta = 1, 2$ and $\delta \geq m+1, m$ respectively, $(\nu_i)_i$ are independent real Brownian motions,  $\lambda_1(0) > \dots > \lambda_m(0) > 0$, $\tau$ is the first collision time  and 
\begin{equation*}
R_0 := \inf\{t, \, \lambda_m (t) = 0\}
\end{equation*}
is the first hitting time of zero. 
Set $r_i := \sqrt{\lambda_i}$, then, for $t <  \tau \wedge R_0$: 
\begin{align*}
dr_i(t) &= d\nu_i(t) + \frac{1}{2r_i(t)}\left[\beta \delta -1 + \beta \sum_{j \neq i}\frac{r_i^2 + r_j^2}{r_i^2-r_j^2}\right] dt \\&
= d\nu_i(t) + \frac{k_0}{r_i(t)}dt + k_1 \sum_{j \neq i}\left[\frac{1}{r_i(t) - r_j(t)} + \frac{1}{r_i(t) + r_j(t)}\right] dt
\end{align*}
with $2k_0 = \beta(\delta - m + 1) -1,\, 2k_1 = \beta$. 
Now consider a $B_m$-type root system 
\begin{equation*}
R = \{\pm e_i,1 \leq i \leq m,\, \pm e_i \pm e_j,\, 1 \leq i < j \leq m \}
\end{equation*}
and take the following positive and simple systems  
\begin{equation*}
R_+ = \{ e_i,1 \leq i \leq m,\, e_i \pm e_j,\, 1 \leq i < j \leq m\}, \, S = \{ e_i - e_{i+1},\, 1 \leq i \leq m-1,\, e_m\}.
\end{equation*}
Then $r = (r_i(t), 1 \leq i \leq m)_{t < \tau \wedge R_0}$ is a radial Dunkl process associated with the root system $R=B_m$ and special multiplicity values. With regard to the exhibited examples and especially to the first one ($R=B_1$), one strongly believes that a similar existence and uniqueness result holds true. This is actually the content of our main result: the radial Dunkl process is the unique strong solution of 
\begin{equation*}
dX_t^W = dB_t  - \nabla \Phi(X_t^W) dt , \quad X_0 \in \overline{C} 
\end{equation*}
where $\Phi(x) = -\sum_{\alpha \in R_+}k(\alpha)\ln(\langle\alpha,x\rangle)$ subject to $k(\alpha) > 0$ for all $\alpha \in R$. This result was proved in \cite{Scha} using Dirichlet forms and a well posed martingale problem. Our proof relies rather on important results due C\'epa and L\'epingle on stochastic differential equations with singular drifts (\cite{Cepa}). It also gives a positive answer to a conjecture announced by Gallardo and Yor in \cite{Chy} stating that for any starting point $X_0 \in V$, the total length of the jumps performed by a Dunkl process up to time $t$ is almost surely finite provided that $k(\alpha) > 0$ for all $\alpha \in R$. Note that since the Dunkl process only jumps along hyperplanes orthogonal to roots, then this conjecture is equivalent to the finiteness of jumps in a given (positive) root direction. Our arguments allow also to derive further results on the mean number of these jumps according to the value of the multiplicity function on each orbit of a given (positive) root (private communication with D. L\'epingle). 

Finally, we give another proof of the following result (see \cite{Chy} for the original proof): the first hitting time $T_0$ of $\partial C$ by $X^W$ 
\begin{equation*}
T_0 := \inf\{t > 0, X_t^W \in \partial C\}
\end{equation*} 
is finite almost surely when $0 \leq k(\alpha) < 1/2$ for at least one simple root $\alpha \in R_+$. We shall even prove the more precise result that if $k(\alpha) < 1/2$ for a simple root $\alpha$, then $X^W$ hits almost surely the facet corresponding to that simple root (part of $\alpha^{\bot}$ lying in $\partial C$). The proof uses ideas from \cite{Cepa} for $A$-type root systems ($T_0 = \tau$ in this setting).
For a better exposition of the paper, the proof of this last result precedes proofs of results concerned with jumps.      

\section{Radial Dunkl Process : Existence and Uniqueness of a strong solution} 
\begin{teo}
\label{exi}
Let $R$ be a reduced root system and recall that:
\begin{equation*}
\Phi(x) = -\sum_{\alpha \in R_+}k(\alpha)\ln(\langle \alpha,x\rangle) := \sum_{\alpha \in R_+}k(\alpha)\theta(\langle\alpha,x\rangle), \quad x \in C, 
\end{equation*}
where $k(\alpha) > 0$ for all $\alpha \in R_+$. Then $X^W$ is the unique strong solution of 
\begin{equation}\label{DE}
dY_t = dB_t - \nabla \Phi(Y_t) dt,  \quad Y_0 \in \overline{C}, \, t \geq 0, 
\end{equation}
where $B$ is a Brownian motion in $V$ and $Y$ is a continuous $\overline{C}$-valued process. 
\end{teo}
\begin{nota}
According to this theorem, the eigenvalues of Wishart and Laguerre processes are unique strong solutions of \eqref{Lag} for any $\lambda_1(0) \geq \lambda_2(0) \geq \dots \geq \lambda_m(0) \geq 0$ provided that $\beta > 0, \delta > m-1 + (1/\beta)$, thereby known results from matrix theory are improved (\cite{Bru}, \cite{Dem}). In this stream, we point the interested reader to the recent paper \cite{GM} where the set of parameters for which the strong uniqueness holds is enlarged, yet the process is not allowed to start at any boundary value. 
\end{nota}

{\it Proof}: From Theorem 2. 2 in \cite{Cepa1}, the SDE: 
\begin{equation}\label{DE1}
dY_t = dB_t - \nabla \Phi(Y_t) dt + n(Y_t)dL_t, \quad Y_0 \in \overline{C} 
\end{equation} 
where $n(x)$ belongs to the set of  unitary inward normal vectors to $C$ at $x \in V$ defined by 
\begin{equation}
\langle x - a, n(x) \rangle \, \leq 0, \quad a \in \overline{C}, 
\label{Ines}
\end{equation}  
and $L$ is the boundary process satisfying: 
\begin{equation*} 
dL_t = {\bf 1}_{\{Y_t \in \partial C\}} dL_t, 
\end{equation*} 
has a unique strong solution for all $t \geq 0$. 
Moreover: 
\begin{eqnarray}\mathbb{E}\left[ \int_0^T {\bf 1}_{\{Y_t \in \partial C\}} dt \right] & = & 0 \label{Slim},\\ 
\mathbb{E}\left[\int_0^T |\nabla \Phi(Y_t)| dt\right] & < & \infty \label{lotfi}
\end{eqnarray}  
for all $T > 0$. All what we need is to prove that $(L_t)_{t \geq 0}$ vanishes. To proceed, we need two Lemmas. 

\begin{lem} 
\label{Hamdi} 
Set $dG_t : = n(Y_t)dL_t$. Then, for all $\alpha \in R_+$,  
\begin{equation*} 
{\bf 1}_{\{\langle Y_t,\alpha\rangle = 0\}} \langle dG_t, \alpha\rangle = 0. 
\end{equation*}
\end{lem}
{\it Proof}: The proof is roughly an extension to arbitrary root systems of the one given in \cite{Cepa} for $R= A_{m-1}$ . In order to convince the reader, we provide an outline. 
On the one hand, the occupation density formula yields: 
\begin{equation*}
\int_0^{\infty}L_t^a(\langle \alpha,Y\rangle)|\theta^{'}|(a)da = \langle \alpha,\alpha\rangle  \int_0^t |\theta^{'}|(\langle\alpha,Y_s\rangle) ds
\end{equation*}
where $L_t^a(\langle \alpha,Y \rangle)$ is the local time up to time $t$ at the level $a \geq 0$ of the real continuous semimartingale $\langle \alpha, Y \rangle \geq 0$ (\cite{Rev}). 
On the other hand, the following inequality holds (instead of (2.5) in \cite{Cepa1}) for all $a \in C$: 
\begin{align*}
\langle \nabla &\Phi(x),x-a \rangle = \sum_{\alpha \in R_+}k(\alpha) \theta^{'}(\langle \alpha,x\rangle)\langle \alpha,x-a\rangle \\& 
\overset{(\star)}{\geq} \sum_{\alpha \in R_+}k(\alpha)[b_{\alpha}|\theta^{'}|(\langle\alpha,x\rangle) - c_{\alpha}\langle\alpha,x-a\rangle - d_{\alpha}] \\&
\geq \min_{\alpha \in R_+}(b_{\alpha}k(\alpha))\sum_{\alpha \in R_+}|\theta^{'}|(\langle\alpha,x\rangle) - |x-a|\sum_{\alpha \in R_+}k(\alpha){c_{\alpha}}|\alpha| - 
\sum_{\alpha \in R_+}k(\alpha)d_{\alpha}
\\& := A \sum_{\alpha \in R_+}|\theta^{'}|(\langle\alpha,x\rangle) - B|x-a| - C
\end{align*}
by Cauchy-Schwarz inequality, where in $(\star)$, we used equation (2.1) in \cite{Cepa1}: let $g$ be a convex $C^1$-function on an open convex set $D \subset \re^m$, then $\forall a \in D$, there exist
$b,c,d > 0$ such that for all $x \in D$: 
\begin{equation*}
\langle \nabla g(x),x-a \rangle  \, \geq \, b|\nabla g(x)| - c |x-a| - d. 
\end{equation*}
 Note also that $A > 0$ since $b_{\alpha}k(\alpha) > 0$ for all $\alpha \in R_+$. Then, the continuity of $Y$ and (\ref{lotfi}) yield: 
\begin{equation*}
\int_0^t |\theta^{'}(\langle \alpha,Y_s \rangle)|ds< \infty
\end{equation*} 
almost surely, which in turn entails:
\begin{equation*}
\int_0^{\infty}L_t^a(\langle \alpha,Y \rangle)|\theta^{'}(a)|da < \infty. 
\end{equation*}
Thus, $L_t^0(\langle \alpha,Y\rangle) = 0$ since the function $a \mapsto |\theta^{'}(a)|$ is not integrable at $0$. 
The next step consists in using Tanaka formula to compute 
\begin{align*}
dZ_t & := d[\langle\alpha,Y_t\rangle - (\langle\alpha,Y_t\rangle)^+]
\\& = {\bf 1}_{\{\langle\alpha,Y_t\rangle = 0\}}\langle\alpha,dB_t \rangle  -  {\bf 1}_{\{\langle \alpha,Y_t\rangle = 0\}} \langle\alpha, \nabla \Phi(Y_t) \rangle dt +   {\bf 1}_{\{\langle 
\alpha,Y_t \rangle = 0\}}\langle \alpha, dG_t \rangle
\end{align*}
for $\alpha \in S$. It is obvious that the second term vanishes. The first vanishes too since it is a continuous local martingale with null bracket (occupation density formula). 
As $Y_t \in \overline{C}$, then $dZ_t = 0$ a.s. which gives the result. $\hfill \blacksquare$

\begin{lem}
\label{Abass}
Let $x \in \partial C$. Then $\langle n(x), \alpha \rangle \neq 0$ for some $\alpha \in S$ such that $\langle x,\alpha \rangle = 0$.
\end{lem}
{\it Proof}: assume that $\langle n(x),\alpha\rangle = 0$ for all $\alpha \in S$ such that $\langle x,\alpha\rangle = 0$. The idea is to find a strictly positive constant $\epsilon$ such that $x - \epsilon n(x) \in \overline{C}$ then use the definition of inward normal vectors (see \eqref{Ines} above) to conclude that $n(x) = 0$. Our assumption implies that $\langle x,\alpha \rangle > 0$ for all $\alpha \in S$ such that $\langle n(x),\alpha \rangle \neq 0$. If such simple roots do not exist, that is, $\langle n(x),\alpha\rangle = 0$ for all $\alpha \in S$, then $x - \epsilon n(x) \in \overline C$ for all $\epsilon > 0$. Otherwise, if $\langle n(x),\alpha \rangle < 0$ for these simple roots, then $x- \epsilon n(x) \in \overline{C}$ for all $\epsilon > 0$. 
Finally, if none of these conditions is satisfied, choose 
\begin{equation*}
0 < \epsilon < \min_{\langle x,\alpha \rangle > 0, \langle n(x),\alpha \rangle  > 0}\frac{\langle x,\alpha \rangle}{\langle n(x),\alpha\rangle},  
\end{equation*}
to see that $x - \epsilon n(x) \in \overline{C}$. Substituting $a = x -\epsilon n(x)$ in \eqref{Ines} for the three alternatives, it then follows that $n(x)$ is the null vector, contradiction. 
$\hfill \blacksquare$\\ 
Now we proceed to end the proof of Theorem \ref{exi}. Lemma \ref{Abass} asserts that 
\begin{equation*}
\{Y_t \in \partial C\} \subset \cup_{\alpha \in S}\{\langle Y_t,\alpha\rangle = 0, \langle n(Y_t),\alpha\rangle \neq 0\}
\end{equation*}
for any time $t$. It follows that 
\begin{align*}
0 &\leq L_t \leq \sum_{\alpha \in S} \int_0^t {\bf 1}_{\{\langle Y_s,\alpha\rangle = 0, \langle n(Y_s),\alpha\rangle \neq 0\}}dL_s
\\& = \sum_{\alpha \in S} \int_0^t \frac{1}{\langle n(Y_s),\alpha\rangle}{\bf 1}_{\{\langle Y_s,\alpha\rangle = 0, \langle n(Y_s),\alpha\rangle \neq 0\}}\langle n(Y_s),\alpha\rangle dL_s = 0
\end{align*}
by Lemma \ref{Hamdi}. $\hfill \blacksquare$ 
\begin{note} 
1/When $m=1$, $(X_t^W)_{t \geq 0}$ is a Bessel process of dimension $\delta = 2k_0 + 1$ and  $k_0 > 0 \Leftrightarrow \delta > 1$. It is a known fact that the local time vanishes (see Ch. XI in \cite{Rev}).\\
2/Since the boundary $\partial C$ of $C$ is a cone, then for $x \in \partial C$, one has $cx \in \partial C$ for any $c \geq 0$. Letting $a = 0$ and $a = cx$ for $c > 1$ in \eqref{Ines}, 
one gets $\langle n(x),x\rangle = 0$ so that $n(x) \in x^{\bot}$ and $\la a,n(x) \ra \geq 0$ for any $a \in \overline{C}$. Moreover,  if $Y_t \in \alpha^{\bot}$ for one and only one $\alpha  \in S$, then $\la n(Y_t),\alpha \ra \neq 0$ by Lemma \ref{Abass}\footnote{One can even take in this case $n(Y_t) = \alpha$.} which together with Lemma \ref{Hamdi} yield 
\begin{equation*}   
{\bf 1}_{\{\langle Y_t,\alpha\rangle  = 0\}} \langle dG_t, \alpha\rangle = {\bf 1}_{\{\langle Y_t,\alpha\rangle = 0\}} \la n(Y_t), \alpha \ra dL_t = 0. 
\end{equation*}
Hence, $(L_t)_{t \geq 0}$ vanishes. 
\end{note}

\section{Finiteness of the first hitting time of the Weyl chamber}
Let $T_0 := \inf\{t > 0, X_t^W \in \partial C\}$ be the first hitting time of the Weyl chamber. In \cite{Cepa}, where $R= A_{m-1}$ and $T_0 = \tau$, authors showed that $T_0 < \infty$ a.s. if $0 \leq k_1 < 1/2$. More generally, it was shown that the same holds if $0 \leq k(\alpha) < 1/2$ for some $\alpha \in S$ (see \cite{Chy} p.169). Below we prove the more precise statement: 
\begin{pro}
\label{Sacha}
Let $\alpha_0 \in S$ and $T_{\alpha_0} := \inf\{t > 0,\, \langle \alpha_0,X_t^W\rangle = 0\}$ such that $T_0 = \inf_{\alpha_0 \in S}T_{\alpha_0}$. If $0 \leq  k(\alpha_0) < 1/2$, then  
$(\langle \alpha_0,X_t^W\rangle)_{t \geq 0}$ hits a.s. $0$. Therefore $T_0 < T_{\alpha_0} < \infty$ a.s.
 \end{pro}
{\it Proof}: Assume $k(\alpha) > 0$ for all $\alpha \in R$ and let $\alpha_0 \in S$. Our scheme is a generalization of the one used in \cite{Cepa}, thus we shall show that the process 
$<\alpha_0,X>$ is a.s. less than or equal to a Bessel process of dimension $2k(\alpha_0) +1 : = 2k_0+1$. The result follows from the fact that $2k(\alpha_0) + 1 < 2$ when $k(\alpha) < 1/2$ so that the Bessel process hits zero a.s.. Using (\ref{DE}), one has for all $t \geq 0$
\begin{align*}
d\langle \alpha_0, X_t^W\rangle &= |\alpha_0| d\gamma_t + \sum_{\alpha \in R_+}k(\alpha)\frac{\langle \alpha,\alpha_0\rangle}{\langle\alpha,X_t^W \rangle} dt \\&
= |\alpha_0|d\gamma_t + k_0\frac{|\alpha_0|^2}{\langle\alpha_0,X_t^W \rangle}dt + \sum_{\alpha \in R_+\setminus \alpha_0}k(\alpha)\frac{\langle \alpha,\alpha_0\rangle}{\langle \alpha,X_t^W \rangle} dt.
\end{align*}
Set 
\begin{equation*}
R =  \cup_{j=1}^p R^j,  
\end{equation*} 
where $R^j,\, 1 \leq j \leq p$ denote the conjugacy classes of $R$ under the $W$-action, then 
\begin{equation*}
R_+ = \cup_{i=1}^p R_+^j \end{equation*}
so that: 
\begin{equation*}   
d\la\alpha_0, X_t^W\ra   = |\alpha_0|d\gamma_t + k_0\frac{|\alpha_0|^2}{\langle \alpha_0,X_t^W \rangle}dt + \sum_{j=0}^p k_j\sum_{\alpha \in R_+^j\setminus \alpha_0}
\frac{\langle\alpha,\alpha_0\rangle}{\langle \alpha,X_t^W\rangle}dt.
\end{equation*}
For a conjugacy class $R^j$ and $\alpha \in R^j$, if $\langle \alpha, \alpha_0\rangle := a(\alpha) > 0$ then, it is easy to check that $\langle \sigma_0(\alpha), \alpha_0\rangle = -a(\alpha)$ where $\sigma_0$ is the reflection with respect to the orthogonal hyperplane $H_{\alpha_0}$. Note that $\sigma_0(\alpha)$ belongs to the same conjugacy class of $\alpha$ and that 
$\sigma_0(\alpha) \in R_+ \setminus \alpha_0$ for $\alpha \in R_+ \setminus \alpha_0$ (see Proposition 1. 4 in \cite{Hum}). Hence, 
\begin{equation*}
d\langle\alpha_0, X_t^W\rangle   = |\alpha_0|d\gamma_t + k_0\frac{|\alpha_0|^2}{\langle\alpha_0,X_t^W\rangle} dt- \sum_{j=0}^p k_j\sum_{\substack{\alpha \in R_+^j\setminus \alpha_0 \\  a(\alpha) >0}}\frac{a(\alpha)\langle \alpha - \sigma_0(\alpha), X_t^W \rangle}{\langle\alpha,X_t^W \rangle\, \langle \sigma_0(\alpha),X_t^W\rangle}dt. 
\end{equation*}
Furthermore, 
\begin{equation*}
\alpha  - \sigma_0(\alpha) = 2\frac{\langle \alpha,\alpha_0\rangle}{\langle \alpha_0, \alpha_0\rangle}\alpha_0 \quad \textrm{so that} \quad \langle\alpha  - \sigma_0(\alpha) , X_t^W \rangle
= 2 a(\alpha)\frac{\langle \alpha_0, X_t^W\rangle}{|\alpha_0|^2}.   
\end{equation*} 
Consequently: 
\begin{equation} \label{scal}
d \langle\alpha_0, X_t^W\rangle  = |\alpha_0|d\gamma_t + k_0\frac{|\alpha_0|^2}{\langle \alpha_0,X^W_t \rangle}dt  + F_t \,dt 
\end{equation} 
where $F_t < 0$ on $\{T_{\alpha_0} = \infty\}$. Using the comparison Theorem in \cite{Kar} (Proposition 2. 18. p. 293 and Exercice 2. 19. p. 294), one claims that 
$\langle \alpha_0,X_t^W \rangle  \, \leq \, Y_{|\alpha_0|^2t}^x$ for all $t \geq 0$ on $\{T_{\alpha_0} = \infty\}$, where $Y^x$ is a Bessel process defined on the same probability space with respect to the same Brownian motion, of dimension $2k_0+1$ and starting at $Y_0 = x \geq \, \langle\alpha_0,X_0^W\rangle > 0$. This is not possible since a Bessel process of dimension $< 2$ hits $0$ a.s.  (\cite{Rev}, Chap. XI).  $\hfill \blacksquare$
\begin{note}
1/When one allows the multiplicity function $k$ to take zero values at some orbits, the SDE \eqref{DE} holds up to $t < T_0$ (Corollary 6.7 p. 169 in\cite{Chy}). Thus our result remains valid under this assumption. 
\\
2/{\bf Open question}: Given two simple roots $\alpha_1,\alpha_2 \in S$ such that 
\begin{equation*}
0 \leq k(\alpha_1) \neq k(\alpha_2) < 1/2,
\end{equation*} 
that is belonging to different orbits, we already know that $T_{\alpha_i} < \infty, i = 1,2$. Is is possible to compare $T_{\alpha_1}$ and $T_{\alpha_2}$? one way to so that is to seek two processes $R_1,R_2$ such that 
\begin{equation*}
\p( \langle \alpha_1,X_t^W\rangle \, < R_1(t) <\, \langle \alpha_2,X_t^W \rangle \,< R_2(t), \, \textrm{for all } t \geq 0) = 1
 \end{equation*} 
and $R_1,R_2$ hits zero almost surely.
\end{note}

\section{Total length and mean number of jumps of a Dunkl process in a given root direction}
\subsection{Total length and a conjecture by Gallardo-Yor}
Recall that the Dunkl process $X$ is a $V$-valued Markov process with jumps whose projection on the orbits space, identified with $\overline{C}$, is $X^W$. Recall also that a jump at time $t$ can only occur in a direction of a root $\alpha$ provided that 
$X_t = \sigma_{\alpha}(X_{t-}) \neq X_{t-}$ so that the  jump's size is given by 
\begin{equation*}
\Delta X_t := X_t - X_{t-} = 2\frac{\la \alpha, X_{t-}\ra}{\la \alpha,\alpha \ra} \alpha.
\end{equation*}   
Then, it is conjectured that for any strictly positive multiplicity function and any $x \in V$ (see \cite{Chy} p.127)
\begin{equation}\label{Dis1}
\sum_{s \leq t} |\Delta X_s|  = \sqrt{2}\sum_{s \leq t} \sum_{\alpha \in R_+} |\la \alpha, X_{s-}\ra| {\bf 1}_{\{X_s = \sigma_{\alpha} X_{s-} \neq X_{s-}\}}\, < \, \infty   
\end{equation}
almost surely, where we assumed without loss of generality that $|\alpha|^2 = 2$. The conjecture was shown in \cite{Chy} to be true for almost all $x \in V$ and uses  tedious computations together with the Markov property. Here, we use Theorem \ref{exi} to give a quick proof to this conjecture: 
\begin{pro}
For any time $t > 0$, any strictly positive multiplicity function and any $x \in V$, the total length of jumps performed by a Dunkl process up to time $t$ is almost surely finite. 
\end{pro}
{\it Proof}: the strategy consists in carrying first the conjecture to the $W$-invariant setting. We start by recalling that after compensating the discontinuous function displayed in \eqref{Dis1} using the L\'evy kernel (\cite{Chy} p.123), it suffices to prove that  
\begin{align*}
\int_0^t ds \sum_{\alpha \in R_+}\frac{k(\alpha)}{|\la \alpha,X_{s} \ra|}
\end{align*}
has finite expectation for any starting point $x \in V$. To proceed, recall that the semi group density of $X$ is given for $(x,y) \in V^2$ by: 
\begin{equation*}
p_t^k(x,y) = \frac{1}{c_kt^{\gamma + n/2}} e^{-(|x|^2 + |y|^2)/(2t)} D_k\left(\frac{x}{\sqrt t}, \frac{y}{\sqrt t}\right) \omega_k(y)
\end{equation*} 
where 
\begin{equation*}
\gamma := \sum_{\alpha \in R_+} k(\alpha),\quad \omega_k(y) := \prod_{\alpha \in R_+}|\la \alpha,y\ra|^{2k(\alpha)}.
\end{equation*}
Recall also that any $y \in V$ is conjugated to one and only one element, say $y'$, belonging to $\overline{C}$. This gives the decomposition 
\begin{equation*}
V = \cup_{w \in W} w\overline{C}.
\end{equation*}
It follows that 
\begin{align*}
\mathbb{E}_x\left[ \sum_{\alpha \in R_+} \frac{k(\alpha)}{|\la \alpha,X_{s} \ra|}\right] &=  \frac{1}{2}\int_V p_s^k(x,y) \sum_{\alpha \in R} \frac{k(\alpha)}{|\la \alpha,y \ra|}dy 
\\& = \frac{1}{2}\sum_{w \in W}\int_{wC} p_s^k(x,y) \sum_{\alpha \in R} \frac{k(\alpha)}{|\la \alpha,y \ra|}dy 
\\& = \frac{1}{2}\int_C\sum_{w \in W} p_s^k(x,wy) \sum_{\alpha \in R} \frac{k(\alpha)}{|\la \alpha,wy \ra|}dy. 
\end{align*}
Now, from the very definition of root systems, one has $wR = R$ yielding 
\begin{equation*}
\sum_{\alpha \in  R}\frac{k(\alpha)}{|\la \alpha,wy \ra|} = \sum_{\alpha \in R}\frac{k(\alpha)}{|\la \alpha,y \ra|} = 2\sum_{\alpha \in R_+}\frac{k(\alpha)}{|\la \alpha,y \ra|} = 
2\sum_{\alpha \in R_+}\frac{k(\alpha)}{\la \alpha,y \ra}
\end{equation*}
for $y \in C$. Besides, if $ x = w_xx', \,x \in V, w_x \in W, x' \in \overline{C}$, then for $y \in C$ 
\begin{align*}
\sum_{w \in W} p_s^k(x,wy) &= \frac{1}{c_ks^{\gamma + n/2}} e^{-(|x|^2 + |y|^2)/(2s)} D_k^W\left(\frac{x}{\sqrt s}, \frac{y}{\sqrt s}\right) \prod_{\alpha \in R_+}|\la \alpha, wy\ra|^{2k(\alpha)}
\\& = \frac{1}{c_ks^{\gamma + n/2}} e^{-(|x|^2 + |y|^2)/(2s)} D_k^W\left(\frac{x'}{\sqrt s}, \frac{y}{\sqrt s}\right)\prod_{\alpha \in R}|\la \alpha, y\ra|^{k(\alpha)} 
\\& = \frac{1}{c_ks^{\gamma + n/2}} e^{-(|x|^2 + |y|^2)/(2s)} D_k^W\left(\frac{x'}{\sqrt s}, \frac{y}{\sqrt s}\right)\prod_{\alpha \in R_+}\la \alpha, y\ra^{2k(\alpha)} 
\end{align*}
where 
\begin{equation*}
D_k^W(x,y) := \sum_{w \in W}D_k(x,wy) = D_k^W(x',y)
\end{equation*}
is the generalized Bessel function (\cite{Chy}). Thus, we showed that

\begin{align*}
\int_0^t ds \mathbb{E}_x\left[ \sum_{\alpha \in R_+} \frac{k(\alpha)}{|\la \alpha,X_{s} \ra|}\right] &=  
\int_0^t ds\,\mathbb{E}_{x'}\left[\sum_{\alpha \in R_+} \frac{k(\alpha)}{\la \alpha,X_{s}^W \ra}\right]  
\\&=\int_0^tds \mathbb{E}_{x'}\left[\sum_{\alpha \in R_+}k(\alpha)|\theta'|\left(\la \alpha,X_s^W \ra\right)\right]
\end{align*} 
(recall that $\theta(u) := - \ln (u), u > 0$). Next, a slight modification of the proof of Lemma \ref{Hamdi}  gives the inequality 
\begin{align*}
A\sum_{\alpha \in R_+} k(\alpha) |\theta'|(\la \alpha,y \ra) \leq \la \nabla \Phi (y), y - a \ra + B|y-a| + D 
\end{align*}
for all $y,a \in C$ and for some constants $A,B,D > 0$. From the definition of $\Phi$, it is easy to see that 
\begin{equation*}
\la \nabla \Phi (y), y \ra = -\sum_{\alpha \in R_+} k(\alpha) = -\gamma. 
 \end{equation*}
Using the Cauchy-Schwartz inequality, it follows that for all $a \in C$
\begin{equation*}
A\sum_{\alpha \in R_+} k(\alpha) \left|\theta'(\la \alpha,y \ra)\right| \leq  (B|a| + D +\gamma) + |a| |\nabla \Phi (y)| + B|y|.
\end{equation*}
Finally, if $k(\alpha) > 0$ for all $\alpha \in R_+$, then \eqref{Slim} holds and we already know that $|X^W|$ is a Bessel process of index $\gamma +n/2-1$  (\cite{Chy}) so that
\begin{equation*}
\int_0^tds \mathbb{E}_{x'}[|X_s^W|]  <\, \infty.
\end{equation*}
The conjecture is proved. $\hfill \blacksquare$   

\begin{nota}
The proof of the above Corollary shows that  
\begin{equation*}
\mathbb{E}_x\left[\frac{1}{\la \alpha, X_1 \ra}\right] < \infty
\end{equation*} 
for all $\alpha \in R_+$. With regard to $p_1^k(x,y)$ and since the singularities of the function $y \mapsto 1/|\la \alpha,y \ra|$ lie only on $\alpha^{\bot}$, then it seems that the Dunkl kernel $y \mapsto D_k(x,y)$ behaves near to $\alpha^{\bot}$ mostly as $1/|\la \alpha,y\ra|$ for fixed $x \in V$. However, we were not able to come rigorously to such a quite important estimation.  
\end{nota}

\subsection{Mean number of jumps}
It was proved in \cite{Chy} that the mean number of jumps is almost surely finite or almost every starting point provided that the multiplicity function is greater or equal than $1/2$. We are 
now going to prove that the conclusion holds true for any starting point which does not belong to any of the hyperplanes orthogonal to the roots. A converse to this result, proved below too, states that if the multiplicity value on some orbit is less than $1/2$ then for any $x \in V$ the mean number of jumps along any root direction in this orbit has infinite expectation. To proceed, set
\[
  \begin{array}{llll}
   N_t^{\alpha} & := & \sum_{s\leq t}{\bf 1}_{\{X_s=\sigma_{\alpha}X_{s-}\neq X_{s-}\}}   & \mbox{for } t\geq 0 \mbox{ and }
            \alpha\in R_+ \;.
  \end{array}
\]
The positive discontinuous functional $ N_t^{\alpha}$ is compensated by
\[
  k(\alpha)\int_0^t\frac{ds}{\langle\alpha,X_s\rangle^2}
\]
with same expectation. Let $k_i$ be the value of $k$ on a conjugacy class $R^i$.
\begin{pro}
If $k_i>1/2$, for any $\alpha\in R_+^i=R^i\cap R_+$ and $x\notin \cup_{\beta\in R^i}H_{\beta}$, we have $\mathbb{E}_x[N_t^{\alpha}]<\infty$.
\end{pro} 
{\it Proof}:  From Ito's formula,
\[
\begin{array}{lll} 
 \ln\langle\alpha,X_t^W\rangle 
& =  &\ln\langle \alpha,x'\rangle + \displaystyle \int_0^t\frac{\langle \alpha,dB_s \rangle}{\langle \alpha,X_s^W\rangle }\\
 && +\displaystyle \sum_{\beta\in R_+} k(\beta) \int_0^t \frac{\langle \alpha,\beta\rangle}{\langle \alpha,X_s^W\rangle \langle \beta,X_s^W\rangle }ds -\frac{1}{2}\int_0^t\frac{|\alpha|^2}{\langle \alpha , X_s^W \rangle^2}ds.
\end{array}
\]
Setting $\tau_{\varepsilon} = \inf{\{s>0: \inf_{\gamma\in R_+^i}\langle \gamma , X_s^W \rangle <\varepsilon\}}$ 
and taking expectations we get
\begin{equation} \label{log}
 \begin{array}{lll} 
 \mathbb{E}_{x'}[\ln\langle\alpha,X_{t\wedge\tau_{\varepsilon}}^W\rangle] &= & \ln\langle \alpha,x'\rangle + 
\displaystyle \sum_{ \beta\neq \alpha} k(\beta) \mathbb{E}_{x'}\left[\int_0^{t\wedge\tau_{\varepsilon}}  \frac{\langle \alpha,\beta\rangle}{\langle \alpha,X_s^W\rangle \langle \beta,X_s^W\rangle }ds\right]\\
    & &  +\displaystyle (k_i-\frac{1}{2})|\alpha|^2\, \mathbb{E}_{x'}\left[\int_0^{t\wedge\tau_{\varepsilon}}\frac{ds}{\langle \alpha , X_s^W \rangle^2}\right].
\end{array}
\end{equation}
From \cite{Hel} we know that for any reduced root system the function 
\begin{equation*}
p(x):=\prod_{\alpha\in R_+}\langle\alpha,x\rangle
\end{equation*}
is $\Delta$-harmonic and therefore for any $x\notin \cup_{\alpha\in R}H_{\alpha}$,
\[
  \sum_{\alpha\neq \beta  \; 
   \alpha, \beta\in R_+}\frac{\langle\alpha,\beta\rangle}{\langle\alpha,x\rangle\langle\beta,x\rangle}=0.
\]
We apply this identity to the root systems $R^i$, $R^j$, $R^i\cup R^j$ for $i\neq j$ and obtain
\[
  \sum_{\alpha\neq \beta  \; 
   \alpha, \beta\in R_+^i}\frac{\langle\alpha,\beta\rangle}{\langle\alpha,x\rangle\langle\beta,x\rangle}=
   \sum_{\alpha\in R_+^i}\sum_{\beta \in R_+^j}\frac{\langle\alpha,\beta\rangle}{\langle\alpha,x\rangle\langle\beta,x\rangle}=
0.
\]
Then, adding up  equations (\ref{log})
with respect to $\alpha\in R^i_+$, we obtain
\[
  \mathbb{E}_{x'}[\sum_{\alpha\in R_+^i}  \ln\langle\alpha,X_{t\wedge\tau_{\varepsilon}}^W\rangle] =
   \sum_{\alpha \in R_+^i} \ln\langle \alpha,x'\rangle 
    + (k_i-\frac{1}{2})\sum_{\alpha\in R_+^i}\, \mathbb{E}_{x'}\left[ \int_0^{t\wedge\tau_{\varepsilon}}
   \frac{ds}{\langle \alpha , X_s^W \rangle^2} \right].
\]
As $\ln\langle \alpha,X^W_u \rangle \leq |\alpha| |X^W_u|$ and $\mathbb{E}_{x'}[\sup_{u\leq t}|X^W_u|]<\infty$, the left hand 
side is bounded above. When $\varepsilon$ goes to zero, $\tau_{\varepsilon}\rightarrow \infty$ a.s. since the facets 
$\partial C\cap(\cup_{\alpha\in R^i}H_{\alpha})$ are not reached by $X^W$. Moreover,  if $x\notin \cup_{\beta\in R^i}H_{\beta}$,
then $\langle \alpha,x'\rangle > 0$ for any $\alpha\in R^i_+$.
Therefore
\[
  \sum_{\alpha\in R_+^i} \, \mathbb{E}_{x'}\left[ \int_0^t
   \frac{ds}{\langle \alpha , X_s^W \rangle^2} \right] <\infty 
\]
and the same argument as in 5.1 shows that  
\[
\sum_{\alpha\in R_+^i}  \mathbb{E}_x[N_t^{\alpha}]= k_i \sum_{\alpha\in R_+^i} \, \mathbb{E}_{x}\left[ \int_0^t
   \frac{ds}{\langle \alpha , X_s\rangle^2} \right] <\infty\;.
\]
$\hfill  \blacksquare$

There is also a converse result.
\begin{pro}
 If $k \leq 1/2$ on some $W$-orbit $R^j$, then for any $x\in V$ and $t>0$
\[
  \sum_{\alpha \in R^j_+} \mathbb{E}_x[N^{\alpha}_t]= \infty.
\] 
\end{pro}
{\it Proof}: Let $k(\alpha_0)\leq 1/2$ for some $\alpha_0 \in R_+$. Setting $U_t:=\langle \alpha_0, X_t^W \rangle^2$ then (\ref{scal}) yields
\[
  dU_t = 2|\alpha_0| \sqrt{U_t} d\gamma_t +(2k(\alpha_0)+1)|\alpha_0|^2 dt +2F_t U_t dt.
\]
Accordingly, if $Z_t$ is a $2$-dimensional squared Bessel process with driving Brownian motion $\gamma$ and starting point $\langle \alpha_0, x'\rangle^2/|\alpha_0|^2$ then the  method  in Theorem 1.4 of \cite{LG}
proves that $U_t\leq |\alpha_0|^2 Z_t$. But a straightforward computation shows that  $\mathbb{E}_{x'}[1/(Z_s)]=\infty$ for any $s>0$.
This finally entails 
\[
  \mathbb{E}_{x'}\left[\frac{1}{\langle \alpha_0,X_s^W\rangle^2} \right]= \infty.
\]
and the conclusion follows. $\hfill \blacksquare$

{\bf Acknowledgments}: Most of results of this paper are announced in {\it C. R. Acad. Sci. Paris, Ser. I 347 (2009) 1125-1128}. The remaining ones are those concerned with the mean number of jumps up to finite time and it is a pleasure to thank D. L\'epingle who accepted to include them here.

\end{document}